
\input gtmacros.tex
\input pictex
\input gtmonout
\volumenumber{1}
\volumeyear{1998}
\volumename{The Epstein birthday schrift}
\pagenumbers{139}{158}
\published{22 October 1998}
\received{27 October 1997}
\papernumber{7}

\reflist
\refkey\B {\bf M Bestvina}, {\bf M Feighn}, {\it Bounding the complexity of
simplicial group actions on trees}, Invent. Math. 103 (1991) 449--469

\refkey\BF {\bf M Bestvina}, {\bf M Feighn}, {\it A counterexample to generalized
accessibility}, from: ``Arboreal Group Theory'',  R\,C Alperin (editor),
MSRI Publications 19,\break Sprin\-ger--Verlag (1991) 133--142 

\refkey\BO {\bf B\,H Bowditch}, {\it Notes on Gromov's hyperbolicity criterion for
path-metric spaces}, from: ``Group Theory from a Geometrical Viewpoint'', E Ghys, A Haefliger and A Verjovsky (editors), World Scientific (1991) 64--167

\refkey\DD {\bf W Dicks}, {\bf M\,J Dunwoody}, {\it Groups acting on graphs},
Cambridge University Press, Cambridge, UK (1989)

\refkey\D {\bf M\,J Dunwoody}, {\it The accessibility of finitely presented
groups}, Invent. Math. 81 (1985) 449--57

\refkey\DU {\bf M\,J Dunwoody}, {\it Groups acting on $\bf R$--trees},
Commun. Algebra 19 (1991) 2125--2136\stdspace {\it Erratum}, 
Commun. Algebra 19 (1991) 3191

\refkey\DUN  {\bf M\,J Dunwoody}, {\it An inaccessible group}, from: 
``The Proceedings of Geometric Group Theory 1991'',
G\,A Niblo and M Roller (editors), LMS Lecture Note Series 181, 
Cambridge University Press (1993) 75--78

\refkey\DUNW  {\bf M\,J Dunwoody}, {\it Inaccessible groups and
protrees}, J. Pure Appl. Algebra 88 (1993) 63--78

\refkey\DUNWO  {\bf M\,J Dunwoody}, {\it Groups acting on protrees}, J.
London Math. Soc. 56 (1997) 125--136

\refkey\GM {\bf  F\,W Gehring}, {\bf G\,J Martin}, {\it Discrete quasiconformal groups
I}, Proc. London Math. Soc. 55 (1987) 331--358

\refkey\L  {\bf G Levitt}, {\bf F Paulin}, {\it Geometric group actions on
trees}, Amer. J. Math. 119 (1997) 83--102

\refkey\LI  {\bf P\,A Linnell}, {\it On accessibility of groups}, J. Pure Appl.
Algebra 30 (1983) 39--46

\refkey\MS  {\bf G\,J Martin}, {\bf R\,K Skora}, {\it Group actions on
the $2$--sphere}, Amer. J. Math. 111 (1989) 387--402

\refkey\S  {\bf Z Sela}, {\it Acylindrical accessibility for groups},
Invent. Math. 129 (1997) 527--565

\refkey\SK {\bf R\,K Skora}, {\it Splittings of surfaces}, Bull. Amer. Math.
Soc. 23 (1990) 85--90

\refkey\ST  {\bf J\,R Stallings}, {\it Foldings of $G$--trees}, 
from: ``Arboreal Group Theory'', R\,C\break Alperin (editor), 
MSRI Publications 19, Springer--Verlag (1991) 355--368

\endreflist

\title{Folding sequences}
\author{M\thinspace J Dunwoody}
\asciiauthors{M J Dunwoody}
\address{Faculty of Math.Studies\\
University of Southampton\\
Southampton, SO9 5NH, UK}
\email{mjd@maths.soton.ac.uk}

\abstract
Bestvina and Feighn showed that a morphism $S \rightarrow T$ between
two simplicial trees that commutes with the action of a group $G$ can
be written as a product of elementary folding operations.  Here a more
general morphism between simplicial trees is considered, which allow
different groups to act on $S$ and $T$.  It is shown that these
morphisms can again be written as a product of elementary operations:
the Bestvina--Feighn folds plus the so-called ``vertex morphisms''.
Applications of this theory are presented.  Limits of infinite folding
sequences are considered.  One application is that a finitely
generated inaccessible group must contain an infinite torsion
subgroup.
\endabstract

\asciiabstract{%
Bestvina and Feighn showed that a morphism S --> T between
two simplicial trees that commutes with the action of a group G can
be written as a product of elementary folding operations.  Here a more
general morphism between simplicial trees is considered, which allow
different groups to act on S and T.  It is shown that these
morphisms can again be written as a product of elementary operations:
the Bestvina-Feighn folds plus the so-called `vertex morphisms'.
Applications of this theory are presented.  Limits of infinite folding
sequences are considered.  One application is that a finitely
generated inaccessible group must contain an infinite torsion
subgroup.}



\maketitle

\cl{\small\it Dedicated to David Epstein on the occasion of his 60th
birthday} 

\section{Introduction}

A morphism $\phi \co S \rightarrow T$ of finite trees can be written as a
product of elementary {\it folds}, in which two edges with a common
vertex are folded together, and an isomorphism.  Bestvina and Feighn [\B]
have given a generalization of this result. The case they consider is
when $S$ and $T$ are (generally infinite) simplicial $G$--trees for
which $G\backslash S$ and
$G\backslash T$ are finite graphs $T$ is minimal, and $G$ and the edge
 stabilizers of $T$ in
$G$ are finitely generated.  The morphism now becomes a product of
equivariant folds and an isomorphism.  In each such fold a whole orbit
of pairs of edges are folded together. Such an operation is easy to
describe in terms of its effect on the quotient graph $G\backslash S$
and the edge and vertex stabilizers of $S$. These are specified in a {\it
graph of groups} determined by a labelling of the edges and vertices of
$G\backslash S$. In this paper a further generalization is given.  We
now allow different groups to act on $S$ and $T$.  Thus $S$ is a
$G$--tree and $T$ is an $H$--tree and a {\it morphism} $\phi \co S \rightarrow
T$ incorporates a homomorphism $\tilde \phi \co G \rightarrow H$, so that
if we regard $T$ as a $G$--tree via $\tilde \phi $ then $\phi $ is a
morphism of $G$--trees.  As well as the basic folding operations of [\B] it
is also necessary to include {\it vertex morphisms} each of which changes just
one vertex label of the corresponding graph of groups. It is possible to
generalize the Bestvina--Feighn result for the case when $\tilde \phi$
restricts to an injective homomorphism on point stabilizers of $S$.
Under similar restrictions to those specified for a $G$--morphism,
$\phi $ is a product of elementary folds, vertex morphisms and an
isomorphism.  A sequence of such operations  is
called a {\it folding sequence}. We can think of each tree in the
sequence as the realization of a combinatorial tree. The folding  and
vertex morphisms correspond to morphisms of the combinatorial trees. If
we interpret our folding sequence as a folding sequence of combinatorial
trees then we also have to allow subdivision operations.  This is because two
different combinatorial trees may have isomorphic realizations as $\bf
R$--trees.  However if this does happen, then the two trees have
isomorphic subdivisions.  
\par Folding sequences are surprisingly useful.  They yield theoretical
results on decompositions of groups and also provide a way of
constructing groups with strange properties.\par
A $G$--tree $S$ is called reduced if for every edge $e \in ES, G_e =
G_{\iota e}$ implies $\iota e, \tau e$ are in the same orbit.  Let $S$
be a reduced $G$--tree in which every edge group is finite.  Let
$\overline S=
G\backslash S$ and let $(G(-), \overline S)$ be the corresponding
graph of groups. Put $$\eta (S) = \sum _{e \in E\overline S}1/|G(e)|.$$
Linnell [\LI] proved that $\eta (S) \leq 2d_G(\omega {\bf Q}G)-1$ where
$d_G(\omega {\bf Q}G)$ is the minimal number of generators of the
augmentation ideal $\omega {\bf Q}G$ as a ${\bf Q}G$--module.  Linnell's
argument uses norms in $W^*$--algebras. Using a
folding sequence argument we show that $\eta (S) \leq d(G)$, the minimal
number of generators of $G$.  If all the edge stabilizers of $S$ are
trivial, then $\eta (S) = |E\overline S|$ and so $
|E\overline S| \leq d(G)$. This is a weak version of the Grushko--Neumann     
Theorem (see [\DD] or [\ST]).  A stronger version of the Grushko--Neumann Theorem
is obtained by a closer examination of the folding sequence. Stallings 
[\ST] has given a proof of this result using this approach.
\par
Let $G$ be a group.  In [\DUNW] and [\DUNWO] I introduced the idea of a 
$G$--{\it protree}. A {\it splitting sequence} of $G$--trees $T_1, T_2, \dots $
is a sequence such that for each $n$
there is a surjective $G$--map $T_n \rightarrow T_{n-1}$ obtained by
contracting finitely many orbits of edges.
A $G$--protree $P$
arises as the inverse limit of this sequence.  As shown in [\DUNWO], if $ET_n$
is countable for all $n$, then $P$ has a
realization as an $\bf R$--tree, on which $G$ acts by isometries.  In
this $\bf R$--tree the set of branch points intersects each segment in a
nowhere dense subset.  A finitely generated group $G$ is said to be {\it
inaccessible} if there is a splitting sequence of reduced $G$--trees as above,
for which all
edge groups are finite and the number of $G$--orbits of $VT_n$ (or
$ET_n$) tends to infinity.  In this case we obtain a $G$--protree $P$ with
infinitely many orbits of edges. \par
 We prove in Section 3 that if $G$ is finitely generated and $P$ is a
$G$--protree
with countably many edges then the realization of $P$ is a direct limit
of a folding sequence of simplicial $\bf R$--trees.  If the $G_n$--tree
$S_n$ is the $n$--th term of the sequence, then there is a surjective
homomorphism $\tilde \rho _n\co G_n \rightarrow G_{n+1}$ and $G$ is the
direct limit of this system of homomorphisms in the category of groups.
This description of $G$ gives information as to the subgroup structure
of $G$.  In particular either $G \cong G_n$ for all sufficiently
large $n$ or $G$ must contain a subgroup which is the union of
a properly ascending chain of finitely generated subgroups each of which
is contained in an edge stabilizer of $P$.  It follows
that an inaccessible group must contain an infinite locally finite
subgroup. If every edge stabilizer of $S_n$ in $G$ is cyclic (not
necessarily finite), then $G$
must contain a non-cyclic subgroup that is locally cyclic. 
It also follows that if $G$ has an infinite splitting sequence then for
any integer $k$ there is an integer
$n$ such that $G$ contains a nontrivial element which fixes an edge
path in $T_n$ of length at least $k$.  This is also implied by Sela's
results on acylindrical accessibility [\S]. 
\par
Infinite folding sequences were used first by Bestvina and Feighn [\BF] to
give an example of a finitely generated group which had an infinite splitting
sequence in which all edge groups are free abelian of rank 2. 
Subsequently [\DUN], [\DUNW], [\DUNWO] I gave a number of examples of 
inaccessible groups all of
which were constructed (essentially) by means of folding sequences.
\par
Martin and Skora [\MS] have obtained some results on accessible convergence
groups acting on $S^2$.  It is not hard to show that an infinite locally
finite group cannot act as a convergence group on $S^2$.  Hence by
Theorem 4.5 a
finitely generated convergence group acting on $S^2$ must be accessible.
Thus the accessibility condition in the results of Martin and Skora can
be removed (or replaced by a finite generation condition). In particular
it follows that if $G \subset {\rm Hom}(S^2)$ is an orientation
preserving convergence group, then there is a simplicial $G$--tree $T$
such that $G\backslash T$ is a finite graph, all edge stabilizers are
finite, and if $v\in VT$, then the ordinary set $O(G_v)$ is simply
connected.

\section{Folding}

We recall and modify some of the terminology of [\DU] or [\SK].

Let $G$ be a group.  A $G$--{\it tree} $T$ is an $\bf R$--tree with $G$ acting
on the left by isometries.  A $G$--tree is {\it minimal} if it has no
proper $G$--subtree.  
\par
Given an $\bf R$--tree $T$ and $x \in T$, define $B_x = \lbrace [x,y] | y \in
T-x\rbrace.$ Define an equivalence relation on $B_x$ by 
\par $[x,y] \sim [x,z] $ if $ [x,y] \cap [x,z] = [x, w]$ for some $w \in T
- x$.
\par A {\it direction at} $x$ is an element of $B_x/\!\sim $.
 There is a bijection between directions at $x$ and the components of
$T-x$.  A {\it point of reflection} $x$ of a $G$--tree $T$ is a point with two
directions for which there exists $g \in G$ which fixes $x$ and
transposes the two directions at $x$.
We say that $x \in T$ is an {\it ordinary point} if there are 
exactly two directions at $x$ but $x$ is not a point of reflection.
A {\it branch
point } is a point $x$ with more than two directions or equivalently 
for which $T-x$ has more than two components.  A {\it vertex} is a point
which is not an ordinary point.
\par
An $\bf R$--tree is $\it simplicial$ if the set of vertices is
discrete.  For each $x \in T$, let $d(v)$ denote the number of
directions at $x$.
\par
A {\it morphism} from a $G$--tree $S$ to a $G$--tree $T$ is a $G$--map
$\phi \co S \rightarrow T$ such that for each segment $[x,y]$ of $S$ there
is a segment $[x,w] \subset [x,y]$ such that $\phi | _{[x,w]} $ is an
isometry.
\par
Alternatively ([6]) $\phi \co S\rightarrow T$ is a morphism if every segment
has a finite subdivision such that $\phi $ restricts to an isometry on
each segment of the subdivision.
\par 
We generalize the notion of morphism to allow different groups to act on
domain and codomain. Thus if $S$ is a $G$--tree and $T$ is an $H$--tree, 
a {\it morphism} $\phi \co S \rightarrow T$ is a homomorphism $\tilde \phi
 \co G \rightarrow H$, 
and a map $\phi \co S
\rightarrow T$ such that if we regard $T$ as a $G$--tree via $\tilde \phi $
then $\phi $ is a morphism when regarded as a
morphism of $G$--trees.  Such morphisms are discussed in unpublished work
of Skora.
\par
A simplicial $\bf R$--tree $T$ can be regarded as the $\it realization$ of a
simplicial complex, which is a (combinatorial) tree. This will also be
denoted $T$.  Thus $VT$ will
correspond to a non-empty closed discrete subset of the $\bf R$--tree
containing all branch points and $ET$ will be the set of closures of
components  of $T - VT$, where VT is such that each element of $ET$ is a
closed segment the endpoints of which are elements of $VT$.
As a combinatorial tree the
vertices of the edge $e$ are denoted $\iota e, \tau e$. When regarded as
a protree the edges of a tree are regarded as directed pairs.  Usually
an edge of a tree is not directed.  

Bestvina and Feighn [\B] have
shown that any morphism of simplicial $G$--trees is a product of subdivisions
and folds (which are described as operations on the corresponding
combinatorial $G$--trees).  Folds are classified according to their effect on the
quotient graph.  The quotient graph $X = G\backslash T$, together with a
labelling by
subgroups of $G$ which are the stabilizers of a lift of a maximal subtree
$X_0$ of $X$, is known as a graph of groups $(G(-), X)$.  See [\DD]
for an account of this theory.  The basic folds of Type I, II and III
are shown below in Figure 1. These are denoted Type IA, IIA and IIIA in [\B].        
Bestvina and Feighn list other basic folds (Type
IB,IIB, IIIB and IIIC ).  But as they remark, each of these is equivalent to a
combination of Type A folds and subdivisions.  
\fig{1}
\beginpicture\small
\setcoordinatesystem units <1.3mm, 1.3mm>
\multiput {$\bullet $} at 0 20  15 20  40 20  55 20  13 27 /
\putrule from 0 20 to 15 20
\setlinear
\plot 0 20 13 27 /
\putrule from 40 20 to 55 20
\multiput {$V$} at -2 18  38 18 /
\put {$E_1$} at 5 25
\put {$E_2$} at 7 18
\put {$X$} at 12 29
\put {$Y$} at 16 18
\put {$\langle E_1, E_2\rangle$} at 47 22
\put {$\langle X,Y\rangle$} at 55 18
\put {type I} at 27 22
\put {$\Rightarrow $} at 27 20
\multiput {$\bullet $} at 0 0  15 0  40 0  55 0 /
\putrule from 0 0 to 15 0
\putrule from 40 0 to 55 0
\multiput {$V$} at -2 -2  38 -2 /
\put {$E$} at 7.5 2
\put {$X$} at 16 -2
\put {$\langle E, g\rangle$} at 47 2
\put {$\langle X,g\rangle$} at 55 -2
\put {type II} at 27 2
\put {$\Rightarrow $} at 27 0
\multiput {$\bullet $} at 0 -20  15 -20  40 -20  55 -20  /
\circulararc 90 degrees from 15 -20 center at 7.5 -27.5
\circulararc 90 degrees from 0 -20 center at  7.5 -12.5
\putrule from 40 -20 to 55 -20
\multiput {$V$} at -2 -22  38 -22 /
\put {$E_1$} at 7.5 -14.5
\put {$E_2$} at 7.5 -25.5
\put {$X$} at  17 -22
\put {$\langle E_1, E_2\rangle$} at 47 -18
\put {$\langle X,g\rangle$} at 55 -22
\put {Type III} at 27 -18
\put {$\Rightarrow $} at 27 -20
\endpicture
\endfig
In [\DUNWO] I introduced vertex morphisms. A {\it vertex morphism} is a
morphism $\theta \co S \rightarrow T$ of simplicial $\bf R$--trees for which
the only change in the corresponding graph of groups is a change in the
label of one of the vertices.  Thus if the label $U$ is changed to $V$
then there is a surjective homomorphism $\theta _U \co U \rightarrow V$
which restricts to the identity map on subgroups which label incident
edges. For vertex morphisms the group $G$ acting on $S$ is different
from the group $H$ acting
on $T$.
We now generalize the Bestvina--Feighn result to allow different groups
to act on domain and codomain.\par
\noindent
\proc{Theorem} Let  $S, T$ be simplicial $\bf R$--trees.  
Let $G$  act by isometries on $S$  and let $H$
act by isometries
on $T$ so that $G\backslash S$  is finite, and all edge
 stabilizers of $T$
 in $H$ are finitely generated. Also $T$  is a minimal
$H$--tree. Let 
$\phi \co S \rightarrow T$  be a morphism, such that the corresponding
homomorphism $\tilde \phi \co G \rightarrow H$ is surjective, and restricts
to an injective map on
each point stabililizer, then $\phi $  can be
written as
a product of basic folds and vertex morphisms.

\prf We adapt the proof of the Proposition in Section 2 of [\B].  
\par 
{\bf Step 0}\stdspace  We show that if $K$ is a finite simplicial subtree of $S$, then
we can factor $\phi $ as $\gamma \beta $ where $\beta $ is a product of
folds and vertex morphisms and $\gamma $ restricted to $\beta (K)$ is an
embedding.  Also $\tilde \gamma $ is injective on all point stabilizers.
If $\phi |K$ is not already an embedding then we can perform a basic fold
 $\phi _1 \co S \rightarrow S_1$ folding together edges $e_1, e_2$ of $S$ so that
 $\phi (e_1) = \phi (e_2)$ and $e_1, e_2$ are
distinct edges of $X$. The basic fold $\phi _1$
produces at most one new edge group and one new vertex group.  The new
edge group is a subgroup of an existing vertex group. 
It follows that $\tilde \phi _1$ restricts to an injective homomorphism on the
 stabilizers
of all except one
orbit of vertices of $S_1$ and on the stabilizers of all edges.  Thus
there is a vertex morphism $\nu _1\co S_1 \rightarrow T_1$ such that $\phi
 \co S \rightarrow T$ factors $\phi = \phi ^{(1)}\nu _1\phi _1$ as a
morphism of $\bf R$--trees (regarding $T$ as an $H$--tree), and also
$\tilde \phi ^{(1)}\co G_1 \rightarrow H$, the homomorphism corresponding
to $\phi^{(1)}$, retricts to an injective homomophism on all point
stabilizers. Note that $\nu _1\phi _1(K)$ has fewer edges than $K$.  We
can therefore proceed by induction on the number of edges of $K$.
\par \noindent
{\bf Step 1}\stdspace  We now claim that we can factor $\phi $ as $\gamma \beta $ so
that $\gamma $ induces a homeomorphism of quotient graphs, $\tilde
\gamma $ is injective on point stabilizers and $\beta $ is a product of
basic folds and vertex morphisms.  This follows exactly as in the
corresponding argument in [1].  The fact that $T$ is a minimal $H$--tree
and $\tilde \phi $ is surjective, together imply that the induced
morphism $G\backslash S \rightarrow H\backslash T$ is a surjective
simplicial map.  One then uses an induction argument based on the number
of edges of $G\backslash S$, using Step 0.
\par \noindent
{\bf Step 2}\stdspace  Since edge stabilizers in $T$ are finitely generated, we can
use the argument of [1] to show that $\phi $ can be factored $\phi =
\gamma  \beta $ as in Step 1 and in addition $\tilde \gamma $ induces an
isomorphism on all edge stabilizers.
\par \noindent
{\bf Step 3}\stdspace  It follows as in [1] that the $\gamma $ obtained in Step 2 is an
isomorphism. 
\qed

We say that in the $G$--tree $S$ an edge $e \in ES$ is {\it compressible} if
 $G_{\iota e} = G_e$ and $\iota e$ and
$\tau e$ lie in different $G$--orbits. We say that $S$ is {\it reduced}
if it has no compressible edges.  For any $G$--finite $G$--tree $S$ there is a 
reduced $G$--tree $S^*$ for which $VS^*$ is a $G$--retract of $S$: $S^*$
is obtained from $S$ by compressing compressible edges. The retraction
is not, in general, uniquely determined.  The retraction is determined
by a {\it compressing forest} $F$ defined as follows: \par
{\bf (1)}\stdspace $F$ is a subgraph of $G\backslash S = \overline S$.
\par
{\bf (2)}\stdspace The edges of $F$ are oriented (given arrows) so that each
vertex $v \in VF$ has at most one arrow pointing away from it.
\par
{\bf (3)}\stdspace If $e \in EF$ then $G(e) = G(\iota e)$, where the arrow on $e$
points from $\iota e$ to $\tau e$.
\par
{\bf (4)}\stdspace $F$ is maximal with respect to properties (1), (2) and (3).
In particular $VF = V\overline S$.

In each component $c$ of $F$ there is exactly one vertex $v_c$ which
 has no arrow
pointing away from it.  The retraction $S \rightarrow S^*$
corresponding to $F$ induces a retraction
 $\rho \co \overline S \rightarrow \overline S^*, \rho (v) =
v_c, v \in c$. 
\par
It is often convenient to work with reduced trees. 
We know that it is possible 
to factorize a morphism of reduced trees as a product of subdivisions,
folds and vertex morphisms.  Unfortunately subdividing a tree always produces
compressible edges.  We introduce some modified folding operations which
allow us to factorize a morphism of reduced trees so that the
intermediate trees obtained are also reduced.  These modified folds are 
shown in Figures 2 ,3 and 4.

Every morphism of $G$--trees is a product of subdivisions and folds of types
I, II and III. Let $\phi \co S \rightarrow T$ be
such a fold.  Given a compressing forest $F$ in $\overline S$, we will
 describe how to construct a
compressing forest $F'$ in $\overline T$ and describe 
the corresponding induced
morphism $\phi ^* \co
S^* \rightarrow T^*$.  Again these are best described by their effect on
the labelled quotient graphs.
\par   
Subdivision induces an isomorphism on the corresponding reduced trees,
since one enlarges the compressing forest to include half the subdivided edge.
Thus a morphism of reduced trees can always be written as a product of
isomorphisms and the morphisms $\phi ^*\co S^* \rightarrow T^*$ induced by
type I, II and III folds.  These are discussed in detail below.
\par

\fig{2}
\beginpicture\small
\setcoordinatesystem units <1.3mm, 1.3mm>
\multiput {$\bullet $} at 0 40  15 40  40 40  55 40  13 47  68 47  /
\putrule from 0 40 to 15 40
\setlinear
\plot 0 40 13 47 /
\plot 55 40  68 47 /
\put {$E$} at 61 45
\put {$Z^*$} at 67 49
\put {$Y^*$} at 55 38 
\putrule from 40 40 to 55 40
\multiput {$V^*$} at -2 38  38 38 /
\put {$E$} at 5 45
\put {$E_2$} at 7 38
\put {$Z^*$} at 12 49
\put {$Y^*$} at 16 38
\put {$ E_2$} at 47 42
\put {Type 1} at 27 42
\put {$\Rightarrow $} at 27 40
\multiput {$\bullet $} at 0 20  15 20  40 20  55 20    70 20 /
\putrule from 0 20 to 15 20
\put {$\langle X, E_2\rangle$} [b] at 47.5 21
\put {$Y$} [b] at 62.5 21
\put {$\langle X,Y\rangle$} [t] at 55 19
\put {$Y^*$} [tl] at 71 19 
\putrule from 40 20 to 70 20
\multiput {$V^*$} at -2 18  38 18 /
\put {$E_2$} [b] at 7 21
\put {$Y^*$} at 16 18
\put {Type 2} at 27 22
\put {$\Rightarrow $} at 27 20
\multiput {$\bullet $} at 0 0  15 0  40 0 /
\putrule from 0 0 to 15 0
\put {$V^*$} at -2 -2 
\put {$E_2$} at 7.5 2
\put {$Y^*$} at 16 -2
\put {Type 3} at 27 2
\put {$\Rightarrow $} at 27 0
\put {$\langle V^*,Y\rangle$} at 40 -2
\put {$\Rightarrow $} at 27 -20
\multiput {$\bullet $} at 0 -20   40 -20  55 -20   /
\circulararc 360 degrees from 0 -20 center at  7.5 -20
\circulararc 360 degrees from 40 -20 center at  47.5 -20
\multiput {$V^*$} [r] at -1 -20  39 -20 /
\put {$E_2$} [r] at 14 -20
\put {$\langle X, E_2\rangle$} [b] at 47.5 -12
\put {$\langle X,Y\rangle$} [l] at 56 -20
\put {$Y$} [b] at 47.5 -27
\put {Type 4} at 27 -18
\put {$\Rightarrow $} at 27 -20
\multiput {$\bullet $} at 0 -40   40 -40 /
\circulararc 360 degrees from 0 -40 center at  7.5 -40
\circulararc 360 degrees from 40 -40 center at  47.5 -40
\put {$V^*$} [r] at -1 -40
\put {$\langle V^*, Y\rangle$} [l] at  40.5 -40 
\put {$E_2$} [r] at 14 -40
\put {$Y$} [l] at 56 -40
\put {Type 5} at 27 -38
\put {$\Rightarrow $} at 27 -40
\put {$\phantom\bullet$} at 0 -47.5
\endpicture
\endfig

\fig{3}
\beginpicture\small
\setcoordinatesystem units <1.3mm, 1.3mm>
\multiput {$\bullet $} at 0 20  15 20  40 20  55 20  13 27  68 27  70 20 /
\putrule from 0 20 to 15 20
\setlinear
\plot 0 20 13 27 /
\plot 55 20  68 27 /
\put {$X$} at 61 25
\put {$Y$} at 64 22
\put {$X^*$} at 67 29
\put {$Y^*$} at 71 18 
\putrule from 40 20 to 70 20
\multiput {$V^*$} at -2 18  38 18 /
\put {$E_1$} at 5 25
\put {$E_2$} at 7 18
\put {$X^*$} at 12 29
\put {$Y^*$} at 16 18
\put {$\langle E_1, E_2\rangle$} at 47 22
\put {$\langle X,Y\rangle$} at 55 18
\put {Type 6} at 27 22
\put {$\Rightarrow $} at 27 20
\multiput {$\bullet $} at 0 0  15 0  40 0  55 0  13 7  70 0 /
\putrule from 0 0 to 15 0
\setlinear
\plot 0 0 13 7 /
\put {$Y^*$} at 71 -2
\putrule from 40 0 to 70 0
\multiput {$V^*$} at -2 -2  38 -2 /
\put {$E_1$} at 5 5
\put {$E_2$} at 7 -2 
\put {$X^*$} at 12 9
\put {$Y^*$} at 16 -2
\put {$Y$} at 62.5 2
\put {$\langle E_1, E_2\rangle$} at 47 2
\put {$\langle X^*,Y\rangle$} at 55 -2
\put {Type 7} at 27 2
\put {$\Rightarrow $} at 27 0
\multiput {$\bullet $} at 0 -20   40 -20  55 -20   /
\circulararc 360 degrees from 0 -20 center at -7.5 -20
\circulararc 360 degrees from 0 -20 center at  7.5 -20
\putrule from 40 -20 to 55 -20
\circulararc 360 degrees from 40 -20 center at 47.5 -20
\multiput {$V^*$} [r] at -1 -20  39 -20 /
\put {$E_1$} [l] at -14 -20
\put {$E_2$} [r] at 14 -20
\put {$\langle E_1, E_2\rangle$} [b] at 47.5 -19.5
\put {$\langle X,Y\rangle$} [l] at 56 -20
\put {$X$} [b] at 47.5 -11.5
\put {$Y$} [t] at 47.5 -28.5
\put {Type 8} at 27 -18
\put {$\Rightarrow $} at 27 -20
\multiput {$\bullet $} at 0 -40     55 -40   /
\circulararc 360 degrees from 0 -40 center at -7.5 -40
\circulararc 360 degrees from 0 -40 center at  7.5 -40
\circulararc 360 degrees from 55 -40 center at 62.5 -40
\circulararc 360 degrees from 40 -40 center at 47.5 -40
\multiput {$V^*$} [r] at -1 -40  54 -40 /
\put {$E_1$} [l] at -14 -40
\put {$E_2$} [r] at 14 -40
\put {$\langle E_1, E_2\rangle$} [r] at 69 -40
\put {$Y$} [l] at 41 -40
\put {Type 9} at 27 -38
\put {$\Rightarrow $} at 27 -40
\put {$\phantom\bullet$} at 0 -47.5
\endpicture
\endfig

\fig{4}
\beginpicture\small
\setcoordinatesystem units <1.3mm, 1.3mm>
\multiput {$\bullet $} at 0 60  15 60  40 60  55 60  70 60 /
\putrule from 0 60 to 15 60
\putrule from 40 60 to 70 60
\multiput {$V^*$} at -2 58  38 58 /
\put {$E$} at 7.5 62
\put {$X^*$} at 16 58
\put {$\langle E, g\rangle$} at 47 62
\put {$\langle X,g\rangle$} at 55 58
\put {Type 10} at 27 62
\put {$\Rightarrow $} at 27 60
\put {$X$} at 62 62
\put {$X^*$} at 71 58
\multiput {$\bullet $} at 0 40  15 40  40 40  55 40  70 40 /
\circulararc 90 degrees from 15 40 center at 7.5 32.5
\circulararc 90 degrees from 0 40 center at  7.5 47.5
\putrule from 40 40 to 70 40
\multiput {$V^*$} at -2 38  38 38 /
\put {$E_1$} at 7.5 45.5
\put {$E_2$} at 7.5 34.5
\multiput {$X^*$} at  17 38  72 38 /
\put {$\langle E_1, E_2\rangle$} at 47 42
\put {$\langle X,g\rangle$} at 55 38
\put {$X$} at 62 42
\put {Type 11} at 27 42
\put {$\Rightarrow $} at 27 40
\multiput {$\bullet $} at 0 20   40 20    /
\circulararc 360 degrees from 0 20 center at -7.5 20
\circulararc 360 degrees from 0 20 center at  7.5 20
\circulararc 360 degrees from 40 20 center at 47.5 20
\put {$V^*$} [r] at -1 20 
\put {$\langle V^*, g\rangle$} [l] at 40.5 20 
\put {$E_1$} [l] at -14 20
\put {$E_2$} [r] at 14 20
\put {$\langle E_1, E_2\rangle$}  [l] at 56 20
\put {Type 12} at 27 22
\put {$\Rightarrow $} at 27 20
\multiput {$\bullet $} at 0 0   40 0  55 0   /
\circulararc 360 degrees from 0 0 center at -7.5 0
\circulararc 360 degrees from 0 0 center at  7.5 0
\circulararc 360 degrees from 40 0 center at 47.5 0
\multiput {$V^*$} [r] at -1 0  39 0 /
\put {$E_1$} [l] at -14 0
\put {$E_2$} [r] at 14 0
\put {$\langle E_1, E_2\rangle$} [b] at 47.5 8
\put {$\langle X,g\rangle$} [l] at 56 0
\put {$X$} [b] at 47.5 -6.5
\put {Type 13} at 27 2
\put {$\Rightarrow $} at 27 0
\multiput {$\bullet $} at 0 -20   40 -20  55 -20 /
\circulararc 360 degrees from 0 -20 center at  7.5 -20
\putrule from 40 -20 to 55 -20
\multiput {$V^*$} [r] at -1 -20  39 -22 /
\put {$E_2$} [r] at 14 -20
\put {$X$} [b] at 47.5 -19
\put {$\langle X,g\rangle$} [t] at 56 -21
\put {Type 14} at 27 -18
\put {$\Rightarrow $} at 27 -20
\put {$\phantom\bullet$} at 0 -27.5
\endpicture
\endfig

We consider the effect of folds on the
quotient graph $\overline S$ and the quotient reduced graph
$\overline S^*$.   
In the subsequent discussion, and in the diagrams of graphs of groups,
the group corresponding to a given edge or vertex is denoted with the
corresponding capital letter, eg the group corresponding to vertex $v$
is $V$ and the group of $e_1$ is $E_1$.  For any vertex $w$, put $\rho
(w) = w^*$, which therefore has the group $W^*$.  Note that if $W= W^*$
then we can change the arrows on $F$ so that $w$ has no arrows pointing
away from it (by reversing all the arrows on the geodesic from $w$ to
$w^*$).  A change like this has
no effect on $\overline S^*$.\par
We now list the different possibilities for the fold $\phi $ and the
resulting induced fold $\phi^* $
\sh{Type I}\par \noindent
$e_1, e_2 \in F$\par
 We choose the new compressing forest $F'$
to contain all $x \in F, x
\not= e_1, e_2$. Also $e_1, e_2$ fold  to form the
 edge $\langle e_1, e_2\rangle$, which is included in $F'$
with an arrow pointing away from pivot vertex $v$ if and only if one of the
edges $e_1, e_2$ has arrow pointing away from $v$.  
It is easy to check that
$F'$ is a compressing forest and $\phi $ induces an
isomorphism on $S^*$, since the folding takes place in a part of the
tree that is compressed both before and after the fold.  
\ppar
\noindent
{\bf $e_1 \in F, e_2 \notin F$ and $v, y$ in different components
of $F$}\par
  Suppose first that the arrow on $e_1$ goes from $x$ to  $v$.  Then
$X = E_1$. After the fold $F'$ is obtained from  $F$ by deleting $e_1$.  
If $X \leq E_2$, then $\phi ^*$  consists of a composite of Type 1 
folds for each edge $e$ which has a vertex $w$ in the same component of
$F$ as $v$
but for which the arrowed path from $w$ to $v^*$ passes through $x$.
It is important to note that in each such Type 1 fold $E \leq E_2$.
Assume then that $X \not\leq E_2$. If $\langle X, E_2\rangle \not= V^*$ then
after folding the new compressing forest
is obtained by omitting the folded edge and also the edge originally
pointing away
from $y$ if $Y \not= Y^*$. Note that $\langle X, E_2\rangle \not= \langle X,Y\rangle$, since $\langle X,
E_2\rangle$ is a subgroup of $V^*$ but  $Y$ is not
contained in $V^*$. 
Such a fold is called a Type 2 fold.  Note that we can assume $E_2
\not= Y$ in a Type 2 fold, since if $E_2 = Y$, then because 
$v, y$ are in different components of $F$ 
we could get a bigger compressing forest by adding $e_2$. If
$Y = Y^*$, then the induced fold is a combination of Type II folds. 
Similarly if $\langle X, E_2\rangle = V^*$  (so that the folded edge must be added to
$F$) and $Y\not= Y^*$, then the induced fold is
a combination of Type II folds.  If $\langle X, E_2\rangle = V^*$ and $Y = Y^*$ then
the induced fold is a Type 3 fold.     
\par
If the arrow on $e_1$ goes from $v$ to $x$, then the fold
produces a compressible edge which can be included in the the new
compressing forest with arrow going from $v$ to $\langle x,y\rangle$.  If there are
arrows in $F$ pointing away from $x$ and $y$ then these edges must be
omitted from the new compressing forest.
If $X \not= X^* (=V^*)$ and $Y \not= Y^*$, the effect on $\overline S^*$ is a
 Type 2 fold (with $\langle X, E_2\rangle = X$). 
Note that $E_2$ is a
proper subgroup of $X$, since otherwise we could add $e_2$ to  $F$ and
get a larger compressing forest in $S^*$.  The
induced fold for $X=X^*$
and $Y\not= Y^*$ is a combination of Type II folds (with $y$ as the
pivot vertex instead of $v$).  The vertex which is initially labelled
$V^*$ finishes with label $\langle V^*, Y\rangle$ and the vertex with label $Y^*$ is
unchanged. The folded edge
becomes a vertex if $X = X^*$ and $Y=Y^*$.  Thus we have a Type 3 fold.
\ppar
{\bf $e_1 \in F, e_2 \notin F$ and $v, y$ in the same component of $F$}\par
We can assume $E_2
\not= Y$, since if $E_2 = Y$ we could change $F$ so that it
included both $e_1$ and $e_2$ which is a case already considered.  To
see this note that $v^* = y^*$.  If there is no edge of 
$F$ pointing away from $y$ then $v^* = y$ and $V= Y$ and we can 
change arrows so that there
is an edge in $F$ pointing away from $y$.  Now change $F$ so that it
includes $e_2$ and omits this edge. Thus we can assume $E_2 \not= Y$.
The analysis for this case is now
similar to that when $v, y$ are in different components.  The induced
fold is of Type 4 if $\langle X, E_2\rangle \not= V^*$ and of Type 5 if $\langle X, E_2\rangle =
V^*$.  Note that, since the part of the graph of groups we are concerned
with in this case is not a tree, it cannot be assumed that all the edge
labels are subgroups of the incident vertex labels.  Thus in a Type 4
fold, $Y$ is not
assumed to be a subgroup of $V^*$---it is conjugate to a subgroup of
$V^*$.  There is no analogous case to Type 3.

\ppar \noindent
{\bf $e_1 \notin F, e_2 \notin F, v, x, y$ in distinct components
of $F$}
\par If either $X = X^*$ or $Y= Y^*$, then we can change the arrows on
$F$ so that either $x$ or $y$ has no edges pointing away from it.
Thus if $F$ contains edges pointing away from both $x$ and $y$, then we
can assume $X \not= X^*$ and $Y \not = Y^*$.  In this case we must
omit at least one of these edges from
$F$ after the fold. 
If $\langle X, Y\rangle \not= X$ then we must omit the edge of $F$
with initial vertex $x$.  Similarly if $\langle X, Y\rangle \not= Y$ then we must omit the
edge of $F$ with initial vertex $y$. If $\langle X,Y\rangle = X = Y$ then we need
only omit one of the two edges, and we can choose either.
First consider the case when both edges are omitted.
The fold in this case is a Type 6 fold if $V^* \not= \langle E_1,
E_2\rangle$.   Note that $E_1 \not= X$ and $E_2 \not= Y$, since otherwise we could add
$e_1$ or $e_2$ to $F$, contradicting its maximality.                    
If $V^* = V = \langle E_1, E_2\rangle$ then the folded edge is compressible and can
be added to $F$.  The induced fold in this case is a combination of Type
II folds: first operating on the edge $e_1$ by increasing $E_1$ to $X$
and $V^*$ to $\langle V^*, X\rangle$ and then operating on the edge $e_2$ by
increasing $E_2$ to $Y$ and $\langle V^*, X\rangle$ to $\langle X, Y\rangle$.
For any edge of $\overline S$ that is not in $F$ which has a vertex
$w$ for which the path from $w$ to $w^*$ passes through $x$ or $y$ it is
necessary to carry out a Type 1 fold in the reduced graph.  Such an
edge, which initially is incident with $x^*$ in $\overline S^*$ 
becomes incident with  $\langle x,y\rangle$ in $\overline T^*$. \par
Consider now the case when only one edge is omitted. This happens for
example if $X= X^*$ and $Y \not= Y^*$
then the induced fold is of Type 7. If $X=X^*$ and $Y=Y^*$ then
the induced fold is just a Type I fold.
If $v, x, y$ are in different components of $F$ then both   $\langle X,Y\rangle \not=
X$ and $\langle X, Y\rangle \not =
Y$,  since $X \leq Y$ implies $x,y$ are in the same component of $F$.
It follows that the edges after the fold cannot be added to $F$.
\ppar \noindent
{\bf $e_1 \notin F, e_2 \notin F, v, x, y$ not in distinct components
of $F$}

\par  This case is similar to the previous case. We can still assume
that $E_1 \not= X$ and $E_2 \not= Y$.  For if say $E_1 = X$, and $v, x$
are in the same component of $F$, then either there is an edge in $F$
pointing away from $x$ or $X=V=V^*$ and there is an edge in $F$ pointing
away from $v$.  We can then change $F$ by removing this edge and
replacing it by $e_1$.  Such a change induces an isomorphism on the
reduced graph.  The fold will now involve an edge of $F$ and has been
considered previously.\par
Suppose $v, x,
y$ are all in the
same component of $F$ so that $V^* = X^* = Y^*$ and $\langle X, Y\rangle \not= V^*,
\langle X, Y\rangle \not= X, \langle X, Y\rangle \not= Y$. The induced fold is of Type 8.
  Again it may be necessary to alter by Type 1 folds the incidence of edges to
vertices in $\overline S^*$.  The similarity with 
the case when $v, x, y$ are in different components of $F$ is because in both
cases $F$ is altered in the same way; by omitting the edges pointing
away from the identified vertex $\langle x,y\rangle$.  It may now be the case that
$\langle X,Y\rangle=X$ say. In this case there would be a compressible edge produced
and so we can add an extra edge to $F$ and the induced fold is of Type 9.

\sh{Type II} \par \noindent
{\bf $e \in F$}\par
In such a fold $V \not= E$ and so
the arrow on $e$ must point from $x$ to $v$.  We can include the
folded edge $\langle e, g\rangle$ in $F'$, with arrow pointing from $\langle x, g\rangle$ to $v$.
\par \noindent
{\bf $e \notin F, v,x$ in different components of $F$}\par

We obtain a Type 10 fold for the case when $X \not= X^*$.
 Type 1 folds in  $\overline S^*$ are necessary
corresponding to any edge of
 $\overline S - F$
joined to a vertex $w$ for which the path from  $w$ to $w^*$ passes through
$x$.  If $X=X^*$ then the
induced fold is just a Type II fold.
\ppar \noindent
{\bf  $e \notin F, v,x$ in the same component of $F$}\par
  This is the same as the previous case
except that the
vertices $v^*$ and $x^*$ are identified before and after the folds.
This gives rise to folds of Type 4 and 5.
\par
\noindent
\sh{Type III}\par \noindent
{\bf  $e_1, e_2 \notin F, v, x$ in different components.}\par
We obtain a Type 11 fold when $X \not= X^*$.  Again Type 1 folds may be
neccessary corresponding to any edge of 
$\overline S - F$ joined to
a vertex $w$ for which the path from $w$ to $w^*$ passes through $x$.
If $X=X^*$ then the
induced fold is just a Type III fold.  
\ppar
\noindent
{\bf  $e_1, e_2 \notin F, v, x$ in the same component of $F$}\par
This produces a Type 12 fold if $X = X^* (=V^*)$, and a Type 13 fold
if $X\not= X^*$.
\ppar
{\bf $e_1 \in F, e_2 \notin F$}\par 
In this case, since
$e_2$ has both its vertices in the same component of $F$ it may be the
case that $E_2 = X$.
We obtain a Type 14 fold.

\ppar
We see then that the induced folds in reduced trees may just be a Type
I, II or III fold, but it may be of a type which creates a new
vertex.  For example a Type 6 fold creates a new vertex.\par
Theorem 2.1 can now be adapted for morphisms between reduced trees.

\proc{Theorem} Let  $S, T$ be simplicial reduced $\bf R$--trees.  
Let $G$ act by isometries on $S$ and let $H$ 
act by isometries
on $T$  so that $G\backslash S$ is finite, and all edge
 stabilizers of $T$
in $H$ are finitely generated. Also $T$ is a minimal
$H$--tree. Let 
$\phi \co S \rightarrow T$ be a morphism, such that the corresponding
homomorphism $\tilde \phi \co G \rightarrow H$ is surjective, and restricts
to an injective map on
each point stabililizer, then $\phi $ can be
written as
a product of folds of Type I, II and III or of Types 1 -- 14 and vertex
morphisms and all the intermediate trees are reduced.\endproc

This result enables us to deduce certain bounds on the complexity of
decompositions of finitely generated groups.

\ppar
Let $S$ be a $G$--tree with finite edge stabilizers. Define
$$ \eta (S) = \sum_{e\in E\overline S^*}1/|G(e)|.$$

\proc{Theorem} Let $G$ be a finitely generated group for
which
$d(G)$
is the minimal number of generators, then $\eta (S) 
\leq d(G)$.
\endproc
\prf
Let $W$ be a free group of rank $d(G)$ and let $X$ be the $W$--tree with
one orbit of vertices on which $W$ acts freely, and for which $\eta
(\overline X^*) = d(G)$.  We regard both $X$ and $S$ as simplicial $\bf
R$--trees. 
A surjective homomorphism $\tilde \alpha \co W \rightarrow G$ induces a morphism
$ \alpha \co X \rightarrow S$.  By Theorem 2.1 $\alpha $ is a
product of basic folds and vertex morphisms. We consider the induced
folds on the reduced trees.
  One can check without too much difficulty
that  $\eta (S)$ does not increase for each of the induced
folds described above. For example, for a fold of Type 6  $$\eta (S)
- \eta (T) = {1\over |E_1|} 
+ {1\over |E_2|} - {1\over |\langle E_1, E_2\rangle|} - {1\over |X|} - {1\over |Y|}.$$
 We can assume $|E_1| \leq |E_2|$. Also we know that $E_1
< X$ and $E_2 < Y$.  Thus ${1\over |X|} \leq {1\over 2|E_1|}$ and
${1\over |Y|} \leq {1\over 2|E_2|} \leq {1\over 2|E_1|}$, so that
${1\over |X|} + {1\over |Y|} \leq {1\over |E_1|}$.  Also ${1\over |\langle E_1,
E_2\rangle|} \leq {1\over |E_2|}$. It is clear in this case
 that  $\eta (S) - \eta
(T) \geq 0.$ Similar arguments show that $\eta (
S)$ does
not increase in each of the other cases. A vertex morphism will
leave edge groups unchanged and cannot increase $\eta (S)$. The theorem
is proved.  
\qed
Let us consider the case when $G$ is a finitely generated group and $S$
is a $G$--tree with trivial edge stabilizers.  In this case $\eta (S)
 = |E\overline S^*|)$, and so we see that the number of
edge orbits in a minimal reduced $G$--tree is bounded by $d(G)$.  In fact
we obtain stronger versions of the Grushko--Neumann Theorem by examining
the folding sequence in this case. Thus we obtain the following theorem,
first obtained in [\DD , I, 10.3].

\proc{Theorem} Let $S$ be a $G$--tree and let $T$
be a reduced
 minimal $H$--tree for
which $G$ acts freely on $ES$ and $H$  acts freely on
$ET$. Also suppose 
$H$ is finitely generated.  Let $\alpha \co S \rightarrow T$ be a
morphism.  
If $\tilde \alpha\co G \rightarrow H$ is surjective then there is a
$G$--tree $S'$ and a morphism $\alpha ' \co S' \rightarrow T$
 that induces
an isomorphism $G\backslash S' \rightarrow H\backslash T$  and $\tilde
\alpha '$ induces a surjective homomorphism $G_v \rightarrow H_{\alpha
'(v)}$ for each vertex $v \in VS'$.
\prf  We can carry out vertex morphisms on $S$ and replace each
vertex stabilizer by its image under $\tilde \alpha $.  We will then
have a $\hat G$--tree $\hat S$ for which there is a morphism $\hat \phi \co
\hat S \rightarrow T$ for which the corresponding homomorphism $\hat G
\rightarrow H$ is injective on all point stabilizers. By Theorem 2.1 $\hat
\phi $ is a product of basic folds, subdivisions and vertex morphisms. 
We consider the induced operations on the corresponding reduced trees.
Since all edge groups are trivial, the only possible induced folds that
can occur are
Type I, III, 1, 3 and 5 (with $E_2 = X = \lbrace 1 \rbrace $).
If we carry out the same sequence of induced folds on $S^*$ (leaving out
all the vertex morphisms), we will obtain the $G$--tree $T'$ with the
required properties.\break\hbox{}\qed

\section{Folding sequences}
A folding sequence $(T_n)$, is a sequence of combinatorial trees $T_n$,
satisfying the following conditions:\par
{\bf (a)}\stdspace $T_n$ is a minimal $G_n$--tree, where $G_n$ is finitely
generated.\par
{\bf (b)}\stdspace $T_{n+1}$ can be obtained from $T_n$ either by subdivision, or
by a 
I, II or III fold followed by a vertex morphism.\par
It is often the case that 
corresponding to a folding sequence $(T_n)$ is a {\it folding sequence of
simplicial} $\bf R$--{\it trees}, in which we replace each tree by a
realization  and the folding operations 
induce morphisms of $\bf R$--trees.  In this case we will risk confusion
by  using $T_n$ to
denote both the tree and its realization as an $\bf R$--tree.
There are examples of folding seqences which cannot be realized in the
above way.  For example if for each $n$, $G_{2n-1}\backslash T_{2n-1}$
 is a tree
with two edges $e_{2n-1}, f_{2n-1}$, and $T_{2n}$ is obtained from
$T_{2n-1}$ by subdividing $e_{2n-1}$ into two edges $e_{2n}$ and
$e_{2n+1}$. Then $T_{2n+1}$ is obtained from $T_{2n}$ by a Type I fold,
in which $e_{2n}$ and $f_{2n-1}$ are folded together to form $f_{2n+1}$.
We call such a folding sequence {\it reducible}.  Thus a folding
sequence is reducible if it satisfies the following condition: \par
There exists  $n$, such that for each $m\geq n$ there is 
 a proper subset $E_m \subset ET_m$ which is invariant under $G_m$  and
such that if the folding operation involves an edge of $E_m$ then the
resulting edges are in $E_{m+1}$.  \par
Thus if the folding operation is subdivision of an edge of $E_m$, then
the resulting edges are all in $E_{m+1}$; and if the operation is a Type
I fold in which one of the edges is in $E_m$, then the resulting edge is
in $E_{m+1}$.  In the the above example the folding sequence is
reducible since the sets $E_{2m} = E_{2m-1} = \lbrace f_{2m-1} \rbrace
$, satisfy the above condition.  A folding sequence is {\it irreducible}
if it is not reducible.  

\proc{Theorem} Let $(T_n)$ be an irreducible folding sequence of
combinatorial trees.  The sequence can be realized as a folding
sequence of morphisms of simplicial {\bf R}--trees in which group actions
are by isometries. 
\prf  For each $n$ it is possible to realize the finite folding
sequence $T_1, T_2, \dots ,$ $T_n$ as a folding sequence of morphisms of
simplicial $\bf R$--trees in which the group actions are by isometries.  
To produce such a realization one has to assign a common length to the
edges in each orbit of edges in such a way that the lengths are
compatible with subdivision and so that Type I and Type III folds take place
between edges of equal length.  To achieve such a realization assign
lengths to the edges of $T_n$ and work backwards, noting that the
lengths of edges of $T_i$ are determined by the lengths of edges of
$T_{i+1}$.
  For each $n = 1, 2, \dots $, let $z_n = (\xi
_{n1}, \xi _{n2}, \xi _{n3}, \dots , \xi _{nk})$ be the length of the
edges $e_1, e_2, \dots e_k$ of $G_1 \backslash T_1$ in such a solution. 
We may assume that for each $n, |z_n| = \sum _{i=1}^k\xi _{ni} =1$. By
compactness for the standard $n-1$--simplex $|\sigma _{n-1}|$,
 the sequence $z_n$ has a convergent subsequence.  Let $w_1 =
(\xi _1, \xi _2, \dots \xi _k)$ be the limit point of a convergent
subsequence.   Note that some of values $\xi _i$ may be zero,
but not all.  We now repeat the above process.  For each term of the convergent
subsequence for $w_1$, we can find a vector corresponding to a solution
for the edges of $G_2\backslash T_2$.  The lengths of these vectors is
bounded, since $|w_1| = 1$.  Again by compactness there is a convergent
subsequence converging to $w_2$ and assigning the coefficients of $w_2$
to  $G_2 \backslash T_2$ will be consistent with assigning the
coefficients of $w_1$ to  the lengths of the edges of $G_1\backslash
T_1$. Note that if an edge has been assigned zero length then when
subdivided the parts have zero length and it can be part of a Type I
fold with another edge of zero length. 
Again repeating this process we can eventually assign lengths to all the
edges of $G_n\backslash T_n$ for every $n$ which are consistent with the
folding process.  If all these lengths are non-zero then we have
realized the folding sequence as a folding sequence of simplicial $\bf
R$--trees. 
If some of the edges have zero  length assigned to them, then it is easy
to see that the folding sequence is reducible.  Thus we take $E_m
\subset ET_m$ to be the set of edges assigned zero length.
\qed
It is easy to construct the {\it limit} of such a folding sequence of $\bf
R$--trees. Let $\theta _n = \rho _n\rho
_{n-1} \dots \rho _1 \co T_1 \rightarrow T_{n+1}$. 
Let $d_n$ be the $\bf R$--tree metric in $T_n$. We define a pseudometric
$d$ in $T_1$ by $d(x, y) = \lim _{n \rightarrow \infty}(d_n(\theta _n(x)),
d_n(\theta _n(y)))$. We put $T = T_1/\sim $, where $x \sim y$ if $d(x,y) = 0.$
Clearly $d$ induces a metric on $T$ and this metric space is called the
limit of the folding sequence.
 
I am grateful to Brian Bowditch for supplying the proof of the following
theorem. 
\proc{Theorem} The limit $T$ of the folding sequence $T_n$ is an
$\bf R$--tree.
\prf
 Let $(S, d)$ be a metric space.  In the terminology of
[\BO ], $d$ is a path metric if given any two points $X, Y \in S$ and
$\epsilon > 0$ there is a rectifiable path joining $X$ and $Y$ of length
at most $d(X,Y) + \epsilon$.  Each $(T_n, d_n)$ satisfies the stronger
condition that any two points $X, Y \in T_n$ are joined by a path of
length $d(X, Y)$.  Since for any $x, y \in T_1, 
(d_n(\theta _n(x)),d_n(\theta _n(y)))$ is a
decreasing sequence, it follows easily that $d$ as defined above 
is a path metric on $T$.  It now follows from [\BO ] Proposition 3.4.2 that $T$
is an $\bf R$--tree if given any four points $X, Y, Z, W$ they can be
partitioned into two sets of two elements, without loss of generality,
$\lbrace \lbrace X, y \rbrace , \lbrace Z, W \rbrace \rbrace $, so that 
$$d(X,Y) + d(Z, W) \leq d(X, Z) + d(Y,W) = d(Y,Z)+d(X,W).$$
Since this condition is satisfied in each $T_n$, it must also be
satisfied in $T$.  Thus $T$ is an $\bf R$--tree.
\qed
If $G$ is the direct limit in the
category of groups of the sequence of homomorphisms $\rho _n \co G_n
\rightarrow G_{n+1}$ then there is an action of $G$ on $T$ via
isometries. Suppose in addition the folding sequence satisfies the following
condition\par 
{\bf (c)}\stdspace Two edges of $T_n$ cannot be folded together if they arose as
subdivided parts of the same edge of $T_m$ for some $m < n$.\par

In this case the natural map  $\phi _n\co T_n \rightarrow
T$ restricts to an isometry on each edge of $T_n$ and it is therefore a
morphism of $\bf R$--trees. It is easy to check that $T$ is the direct
limit of the sequence of folding morphisms in the category $\cal T$ of $\bf
R$--trees and morphisms.
\par

As noted above, it is best to describe folding operations in terms of
their effect on the quotient graphs $G_n\backslash T_n$.
Note that (c) applies to $T_n$ and not to $G_n\backslash T_n$. Thus it
is possible for the $n$-th fold in the folding sequence to fold together
edges that arose as 
subdivided edges of $G_m\backslash T_m$ for some $m<n$.  An example of
this is given in
[\DUNW ]. What happens is that, in $T_n$, the edges folded together occur as
subdivided parts of different edges in the same $G_m$--orbit in $T_m$.

\ppar

Let $G$ be a finitely generated group.
Suppose we have an infinite folding sequence with limit $T$ and 
suppose that $\tilde \phi _n \co G_n \rightarrow G$ is not an isomorphism
for any $n$.  This means that the folding sequence must have infinitely
many vertex morphisms.
There is then an induced folding sequence of reduced trees. 
We examine the
induced folds listed above.  For induced folds of type I, III and 
3 there is a decrease in the number of orbits of edges. For a fold of type
12, 13 or 14  there is a decrease in the rank of $H_1(\overline S^*)$ and for a
fold of type 1 there is no change in vertex groups. 
 Thus the sequence must contain infinitely many
induced folds of types other than I, III, 1, 3, 12, 13 or 14.  However each
such induced fold, which is not an isomorphism, produces a new edge
group that properly contains 
one of the old edge groups.  In the situation when the maps $\phi _n\co
T_n \rightarrow T$ 
are morphisms of $\bf R$--trees, for example if condition (c) is
satisfied, each edge stabilizer of $T_n$ fixes an
arc of $T$. Since each $T_n$ has finitely many orbits of edges, using
a graph theoretic argument (K\" onig's Lemma) it is possible to find 
a sequence of edge stabilizers from a subsequence of the $T_n$'s
for which the inclusions are proper.
It follows that $G$
contains a subgroup $H$ that is not finitely generated but every
finitely generated subgroup of $H$ fixes an arc of $T$.  Thus we have
the following result.

\proc{Theorem}  Let the $G$--tree $T$ be the direct limit in
$\cal T$ of the folding
sequence $T_n$ of simplicial trees, where $T$ is a $G_n$--tree.  Then
either there exists $m$ such that $G = G_n$ for all $n \geq m$ or $G$
contains a subgroup $H$ that is not finitely generated but every
finitely generated subgroup of $H$ fixes an arc of $T$.\endproc

In [\DUNW ] I introduced the concept of a $G$--protree.
 Protrees arise naturally in studying
inaccessible groups.  Let $G$ be a finitely generated group.  Let ${\cal
B}(G)$ denote the Boolean ring consisting of all subsets $a \subset G$
of almost invariant sets. Thus $a \in {\cal B}(G)$ if and only if the
sets $a$ and $ag$ are almost equal for every $g \in G$.  In [4] it is
shown that there is a {\it nested} $G$--set $E$ which generates ${\cal
B}(G)$ as a Boolean ring.  The group $G$ is accessible  if and only if
$E$ can be chosen to be $G$--finite, in which case $E$ can be regarded as
the edge set of a simplicial $G$--tree.  If $G$ is inaccessible then $E$
is not $G$--finite.  In this case $E$ is a combinatorial object called a
nice $G$--protree, which has a realization (also called a $G$--protree) as 
an $\bf R$--tree in which the set of branch points intersects each
segment in a nowhere dense subset.  

If $G$ is finitely generated, then
any $G$--tree $T$ is a strong limit of a sequence  $T_n$ of $\bf
R$--trees, where $T_n$ is a $G_n$--tree and the action is geometric, ie it
arises from a foliation on a finite $2$--complex.  See [\L ] for a precise
definition and a proof of the above statement.  However in a geometric
action an orbit which is nowhere dense must be discrete (see [\L ]).  Thus
if $G$ is finitely generated and $T$ is a $G$--protree, then $T$ is a
strong limit of a folding sequence of simplicial trees.  
This gives the following result.
\proc{Theorem} Let $G$ be a finitely generated group and
let $P$ be a nice $G$--protree. Then either \par
{\rm (i)}\stdspace there is
a reduced $G$--tree $T$ such
that for every $v\in VT, G_v$ is finitely generated and fixes a
vertex of $P$ and for
every $e \in ET, G_e$ is finitely generated and fixes an edge of $P$,
\par \noindent or \par
{\rm (ii)}\stdspace the group $G$ 
contains a subgroup $H$ that is not finitely generated but every
finitely generated subgroup of $H$ fixes an edge of  $P$.       
\endproc

Note that if $G$ is finitely presented then $\tilde \phi _n$ must be an
isomorphism for $n$ large and so (i) must hold.  This can be used to
give a proof that
finitely presented groups are accessible. 
This was first proved in [\D ].  We have seen that if $G$ is finitely
generated then we can construct a $G$--protree $P$ corresponding to a nested
set of generators of ${\cal B}(G)$.  There is then a folding sequence
which has limit $P$.  If the situation (i) of Theorem 3.4 prevails then
for each $v \in VT, G_v$ will have at most one end and so $G$ will be
accessible.  Thus if $G$ is inaccessible then condition (ii) must be 
satisfied giving the following result.  
\proc{Theorem} Let $G$ be a finitely generated inaccessible
group. Then $G$ contains an infinite locally finite subgroup.
\prf This follows immediately from Theorem 3.4.
\qed
\proc{Corollary} Let $G$ be a finitely generated discrete
convergence group acting on $S^2$.  Then $G$ is accessible.
\prf By Theorem 3.5 it suffices to show that a locally finite
discrete convergence group must be finite.  Suppose that $H$ is an
infinite locally finite discrete convergence group acting on $S^2$.  By
[\GM ] Theorem 5.11, $L(H)$ (the set of limit points of $H$)
consists of exactly one point $x_0$, which is
fixed by $H$. A finite group of homeomorphisms with a fixed point is 
conjugate in $Hom(S^2)$ to a cyclic or dihedral group acting linearly on
$S^2$ . An increasing
chain of such groups would have to have two fixed points, contradicting
the statement above that there is a unique fixed point.   

\references

\bigskip
{\parskip0pt \small\sl\nl
Faculty of Math.Studies \nl
University of Southampton,\nl
Southampton, SO9 5NH, UK  

\medskip\rm Email:\stdspace\tt
\theemail\par}
\recd

\end